\documentclass[oneside,english]{amsart}
\usepackage[T1]{fontenc}
\usepackage[latin9]{inputenc}
\usepackage{color}
\usepackage{babel}
\usepackage{float}
\usepackage{amstext}
\usepackage{amsthm}
\usepackage{graphicx}
\usepackage[unicode=true,pdfusetitle,
 bookmarks=true,bookmarksnumbered=false,bookmarksopen=false,
 breaklinks=false,pdfborder={0 0 1},backref=false,colorlinks=false]
 {hyperref}
\usepackage{breakurl}

\makeatletter

\floatstyle{ruled}
\newfloat{algorithm}{tbp}{loa}
\providecommand{\algorithmname}{Algorithm}
\floatname{algorithm}{\protect\algorithmname}

\numberwithin{equation}{section}
\numberwithin{figure}{section}
\theoremstyle{plain}
\newtheorem{thm}{\protect\theoremname}
  \theoremstyle{definition}
  \newtheorem{defn}[thm]{\protect\definitionname}
  \theoremstyle{plain}
  \newtheorem{lem}[thm]{\protect\lemmaname}

\usepackage{babel}
\usepackage{algorithm,algpseudocode}

\makeatother

  \providecommand{\definitionname}{Definition}
  \providecommand{\lemmaname}{Lemma}
\providecommand{\theoremname}{Theorem}

\begin{document}

\title{Layout of random circulant graphs}

\author{Sebastian Richter}

\address{Fakultät für Mathematik, Technische Universität Chemnitz, D-09107
Chemnitz, Germany. }

\thanks{\textit{Richter }thanks the Institute of Computer Science of The
Czech Academy of Sciences, and the Czech Science Foundation, under
the grant number GJ16-07822Y, for travel support.}

\email{sebastian.richter@mathematik.tu-chemnitz.de}

\author{Israel Rocha}

\address{The Czech Academy of Sciences, Institute of Computer Science, Pod
Vodárenskou v\v{e}ží 2, 182~07 Prague, Czech Republic. With institutional
support RVO:67985807.}

\thanks{\emph{Rocha} was supported by the Czech Science Foundation, grant
number GJ16-07822Y}

\email{rocha@cs.cas.cz}
\begin{abstract}
A circulant graph $H$ is defined on the set of vertices $V=\left\{ 1,\ldots,n\right\} $
and edges $E=\left\{ \left(i,j\right):\left|i-j\right|\equiv s\left(\textrm{mod}n\right),s\in S\right\} ,$
where $S\subseteq\left\{ 1,\ldots,\lceil\frac{n-1}{2}\rceil\right\} .$
A random circulant graph results from deleting edges of $H$ with
probability $1-p$. We provide a polynomial time algorithm that approximates
the solution to the minimum linear arrangement problem for random
circulant graphs. We then bound the error of the approximation with
high probability.
\end{abstract}



\maketitle

\section{Introduction}

A layout on the graph $G=(V,E)$ is a bijection function $f:V\rightarrow\left\{ 1,\dots,\left|V\right|\right\} $.
Layout problems can be used to formulate several well-known optimization
problems on graphs. Also known as linear ordering problems or linear
arrangement problems, they consist on the minimization of specific
metrics. Such metrics would provide the solution to problems as linear
arrangement, bandwidth, modified cut, cut width, sum cut, vertex separation
and edge separation. All these problems are NP-hard in the general
case.

The Minimum Linear Arrangement (MinLA) problem is to find a function
$f$ that minimizes the sum $\sum_{uv\in E}\left|f(u)-f(v)\right|.$
A layout is also called a labeling, a ordering, or a linear arrangement.
The MinLA is one of the most important graph layout problems and was
introduced in 1964 by Harper to develop error-correcting codes with
minimal average absolute errors. In fact, MinLA appear in a vast domain
of problems: VLSI circuit design, network reliability, topology awareness
of overlay networks, single machine job scheduling, numerical analysis,
computational biology, information retrieval, automatic graph drawing,
etc. For instance, layout problems appear in the reconstruction of
DNA sequences \cite{Karp}, using overlaps of genes between fragments.
Also, MinLA has been used in brain cortex modelling \cite{Mitchison}.
In \cite{DiazPenrose} it is presented a good survey on graph layout
problems and its applications.

The main contribution of this paper is a polynomial time algorithm
that approximates the solution to the MinLA problem for a random circulant
graph. First, a circulant graph $H$ is defined on the set of vertices
$V=\left\{ 1,\ldots,n\right\} $ and edges $E=\left\{ \left(i,j\right):\left|i-j\right|\equiv s\left(\textrm{mod}n\right),s\in S\right\} ,$
where $S\subseteq\left\{ 1,\ldots,\lceil\frac{n-1}{2}\rceil\right\} .$
A random circulant graph results from deleting edges of $H$ with
probability $1-p$. Noticeable, circulant graphs and its random instances
carry a nice shape. The MinLA problem for these graphs is the reconstruction
of that shape.

\begin{figure}[H]
\includegraphics[bb=60bp 550bp 500bp 792bp,scale=0.8]{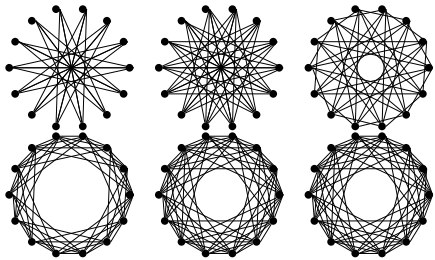}

\caption{Circulant graphs}
\end{figure}

The MinLA problem for circulant graphs is solved in \cite{Rajasingh},
where the authors address the problem of finding an embedding of $G$
into a graph $H$. In that case, $G$ is a circulant graph and $H$
is a cycle graph. Certain circulant graphs are of particular interest.
In \cite{Habibi} it is presented a polynomial time algorithm solving
MinLA of Chord graphs, which is a particular case of circulant graphs.
The main motivation of \cite{Habibi} is an application to topology
awareness of peer-to-peer overlay networks. The solution of \cite{Habibi}
assumes that the Chord graph is complete. However, in real overlay
networks nodes can disconnect at any moment, so the remaining network
can be regarded as a random circulant graph. Thus, the solution to
MinLA for random circulant graphs suits well such applications.

Nevertheless, layout problems for random graphs are significantly
more complicated and usually the solution is an approximation of the
solution of the model graph. The paper \cite{DiazPenrose} is concerned
with the approximability of several layout problems on families of
random geometric graphs. It is proven that some of these problems
are still NP-complete even for deterministic geometric graphs. The
authors present heuristics that turn out to be constant approximation
algorithms for layout problems on random geometric graphs, almost
surely. The authors of \cite{DiazPenrose} remark that their algorithms
use the node coordinates in order to build a layout. That is another
feature we do not require in our problem. Even tough, the random graph
follows a geometric graph model (the circulant structure), we do not
have the coordinates of the random graph in advance. The input random
graph consists of a set of vertices and edges only and we have to
retrieve the circulant layout from that.

Eigenvectors of random matrices are the main tool we use to construct
the layout in our problem. We introduce this idea in \cite{LinearPaper},
where one eigenvector would suffice to recover the structure of a
random linear graph. Here, as we will see, one eigenvector alone is
not enough to encode the whole layout. Fortunately, we can combine
two special eigenvectors to find the linear arrangement. Even tough,
the use of eigenvectors in the same fashion is a common feature of
both methods, here we require some additional technical details that
were not present in \cite{LinearPaper}. Due to the use of angles
between subspaces and SVD decomposition, the technique we use here
differs significantly from \cite{LinearPaper}. There, we pointed
out the generality of such method and here it turns out we need a
more careful analysis. Nevertheless, we have evidence that these methods
can be used to implement a general framework for which layout problems
can be solved in a broader class of random geometric graphs.

The rest of the paper is organized as follow. In section \ref{sec:Main-results}
we define the model matrix, state the algorithm, and the main theorems.
In section \ref{sec:SVD-and-angles} we describe basic properties
of angle between subspaces. Finally, in section \ref{sec:Bounds-and-proofs}
we provide the proofs for the results.

\section{\label{sec:Main-results}Main results}

A circulant matrix $A$ is a matrix that can be completely specified
by only one vector $a$, that appears in the first column of $A$.
The remaining columns are cyclic permutations of $a$ with offset
equal to the column index, i.e., the matrix $A$ is of the following
form 
\[
A=\begin{bmatrix}a_{1} & a_{2} & a_{3} & \dots & a_{n}\\
a_{n} & a_{1} & a_{2} & \dots & a_{n-1}\\
\vdots &  & \ddots &  & \vdots\\
a_{2} & a_{3} & a_{4} & \dots & a_{1}
\end{bmatrix}.
\]
A circulant graph is a graph with circulant adjacency matrix. Let
$H=(V,E_H)$ be a circulant graph with vertex set $V=\{v_{1},\dots,v_{n}\}$
and adjacency matrix $A$, where $[a_{1},...,a_{n}]$ corresponds
to the first row of $A$. We define the set of indices of non-zero
elements in the first half of the row of $A$ as 
\[
N:=\{k:a_{k}=1,1\le k\le\lceil\frac{n-1}{2}\rceil\}.
\]
Equivalently, a circulant graph can be defined as the Cayley graph
of a finite cyclic group.

\begin{figure}[H]
\begin{centering}
\includegraphics[scale=0.15]{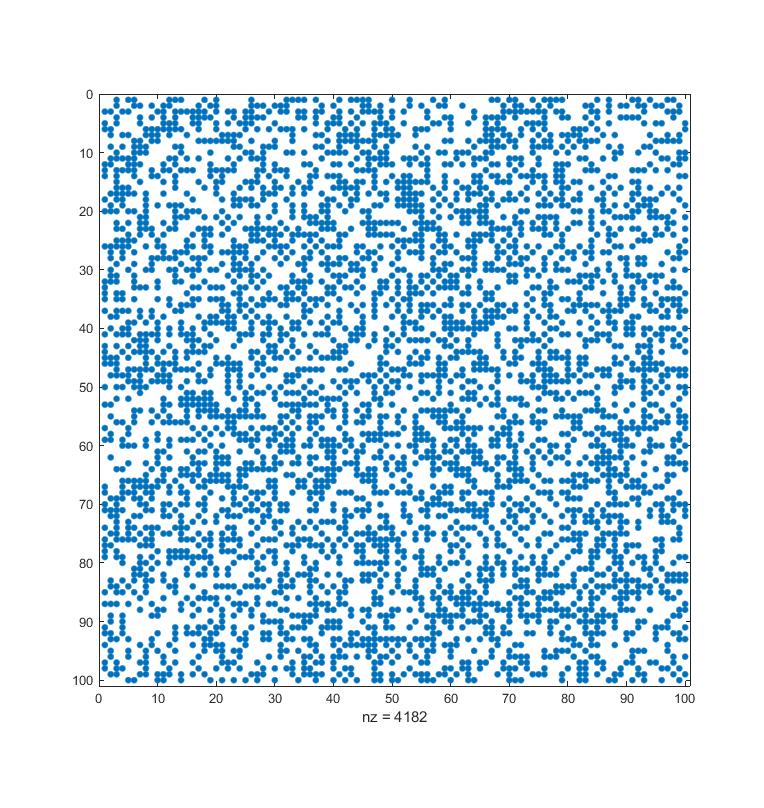}\includegraphics[scale=0.15]{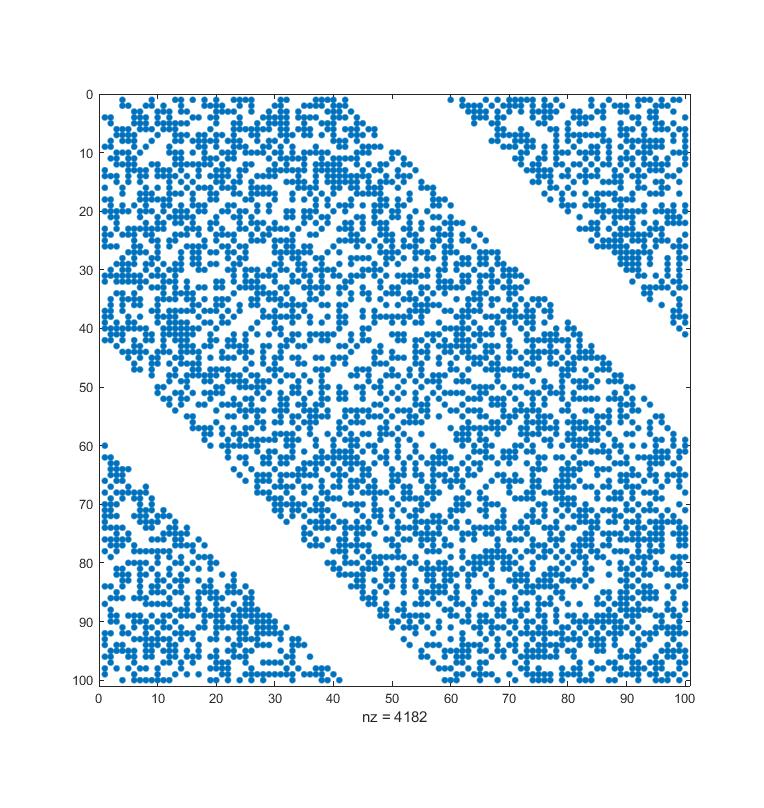}
\par\end{centering}
\caption{Data before and after the correct permutation.}

\label{fig:noisy}

\end{figure}

In this paper, $H$ is referred as the model graph. The random graph
we consider is denoted by $G=(V,E)$ which results from deleting edges
of $H$ with probability $1-p.$ The model matrix $M$ is a circulant
matrix that describes the structure of $H$, where $M=pA$.\begin{figure}[h]
 \begin{minipage}{0.3\textwidth}

\includegraphics[scale=0.25]{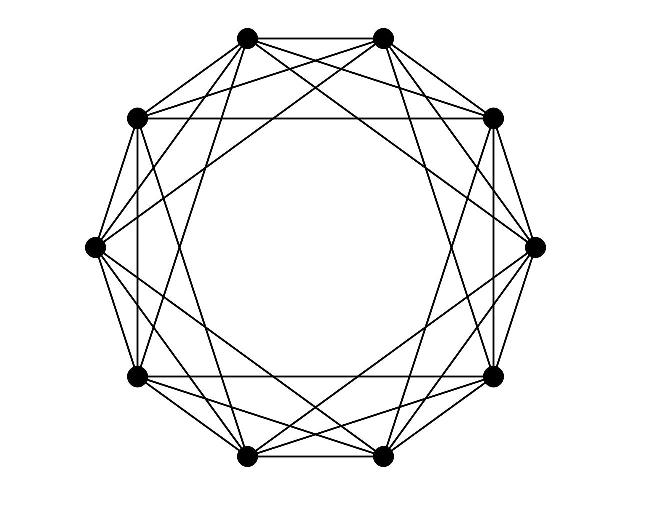}
 \end{minipage}
 \hfill
 \begin{minipage}{0.5\textwidth}{\scriptsize{}$M=p\begin{bmatrix}0 & {\color{blue}1} & {\color{blue}1} & {\color{blue}1} & 0 & 0 & 0 & {\color{blue}1} & {\color{blue}1} & {\color{blue}1}\\
{\color{blue}1} & 0 & {\color{blue}1} & {\color{blue}1} & {\color{blue}1} & 0 & 0 & 0 & {\color{blue}1} & {\color{blue}1}\\
{\color{blue}1} & {\color{blue}1} & 0 & {\color{blue}1} & {\color{blue}1} & {\color{blue}1} & 0 & 0 & 0 & {\color{blue}1}\\
{\color{blue}1} & {\color{blue}1} & {\color{blue}1} & 0 & {\color{blue}1} & {\color{blue}1} & {\color{blue}1} & 0 & 0 & 0\\
0 & {\color{blue}1} & {\color{blue}1} & {\color{blue}1} & 0 & {\color{blue}1} & {\color{blue}1} & {\color{blue}1} & 0 & 0\\
0 & 0 & {\color{blue}1} & {\color{blue}1} & {\color{blue}1} & 0 & {\color{blue}1} & {\color{blue}1} & {\color{blue}1} & 0\\
0 & 0 & 0 & {\color{blue}1} & {\color{blue}1} & {\color{blue}1} & 0 & {\color{blue}1} & {\color{blue}1} & {\color{blue}1}\\
{\color{blue}1} & 0 & 0 & 0 & {\color{blue}1} & {\color{blue}1} & {\color{blue}1} & 0 & {\color{blue}1} & {\color{blue}1}\\
{\color{blue}1} & {\color{blue}1} & 0 & 0 & 0 & {\color{blue}1} & {\color{blue}1} & {\color{blue}1} & 0 & {\color{blue}1}\\
{\color{blue}1} & {\color{blue}1} & {\color{blue}1} & 0 & 0 & 0 & {\color{blue}1} & {\color{blue}1} & {\color{blue}1} & 0
\end{bmatrix}$}\end{minipage}
\caption{Model graph and its model matrix}
\label{fig:model}
\end{figure}\\
Furthermore, let $\hat{M}$ be the adjacency matrix of the random
graph $G$. The entries of $\hat{M}$ correspond to independent Bernoulli
variables, where $\mathbf{P}(\hat{m}_{ij}=1)=m_{ij}.$\begin{figure}[h]
 \begin{minipage}{0.3\textwidth}\includegraphics[scale=0.25]{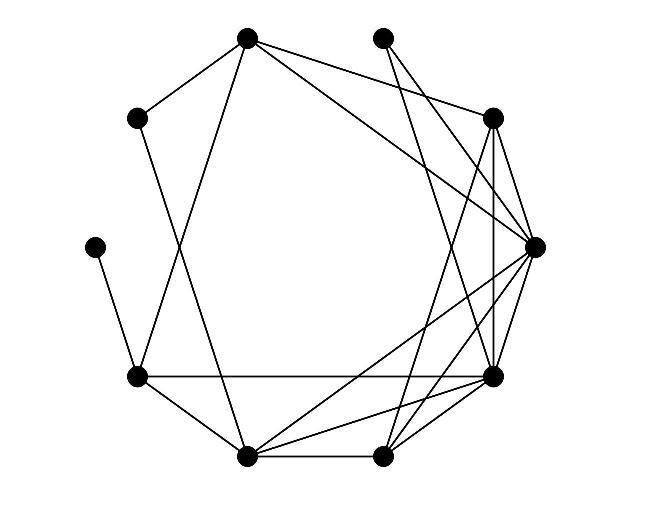}\end{minipage}
 \hfill
 \begin{minipage}{0.5\textwidth}{\scriptsize{}$\hat{M}=\begin{bmatrix}0 & 0 & {\color{blue}1} & {\color{blue}1} & 0 & 0 & 0 & {\color{blue}1} & 0 & {\color{blue}1}\\
0 & 0 & 0 & {\color{blue}1} & {\color{blue}1} & 0 & 0 & 0 & 0 & 0\\
{\color{blue}1} & 0 & 0 & {\color{blue}1} & {\color{blue}1} & {\color{blue}1} & 0 & 0 & 0 & 0\\
{\color{blue}1} & {\color{blue}1} & {\color{blue}1} & 0 & {\color{blue}1} & {\color{blue}1} & {\color{blue}1} & 0 & 0 & 0\\
0 & {\color{blue}1} & {\color{blue}1} & {\color{blue}1} & 0 & {\color{blue}1} & {\color{blue}1} & {\color{blue}1} & 0 & 0\\
0 & 0 & {\color{blue}1} & {\color{blue}{\color{blue}1}} & {\color{blue}1} & 0 & {\color{blue}1} & 0 & 0 & 0\\
0 & 0 & 0 & {\color{blue}1} & {\color{blue}1} & {\color{blue}1} & 0 & {\color{blue}1} & 0 & {\color{blue}1}\\
{\color{blue}1} & 0 & 0 & 0 & {\color{blue}1} & 0 & {\color{blue}1} & 0 & {\color{blue}1} & 0\\
0 & 0 & 0 & 0 & 0 & 0 & 0 & {\color{blue}1} & 0 & 0\\
{\color{blue}1} & 0 & 0 & 0 & 0 & 0 & 0 & {\color{blue}1} & 0 & 0
\end{bmatrix}$}\end{minipage}
\caption{Random graph and its random matrix}
\label{fig:random}
\end{figure}\\
By construction, the labels of vertices in the random graph $G$ corresponds
to the same labels as the graph model in the first figure. However,
in the real world we do not have the labels in advance. We are talking
about a large amount of disorganized data with additional noise. We
only know that this data encodes a circulant structure which is hidden
from us. In such situations, finding the labels for the random graph
can be rather challenging. See Figure \ref{fig:noisy}.

That is precisely the problem we address in this paper: given a graph
that follows a circulant model, find its circular embedding, or rather,
retrieve the correct order of the vertices. We present an algorithm
that solves this problem by using eigenvectors corresponding the second
and third largest eigenvalues $\hat{M}$. The algorithm can be described
as follow.
\begin{algorithm}[H]
\begin{algorithmic}[1] \Require{Random matrix $\widehat{M}$} \\
Compute $\widehat{x}$ and $\widehat{y}$, the eigenvectors for $\lambda_{2}(\widehat{M})$ and $\lambda_{3}(\widehat{M})$\\
Compute the angular coordinate $\varphi_i$ for the point of coordinates $(\widehat{x_i},\widehat{y_i})$\\
Define a permutation $\sigma$ such that $\sigma(i)>\sigma(j)$ iff $\varphi_i\geq\varphi_j$\\
\Return{$\sigma$} \end{algorithmic}

\caption{}
\label{algo:main}
\end{algorithm}
This simple algorithm is shown to return the correct labels with a
bounded error. We quantify the error in terms of a rank correlation
coefficient we introduce. Before, let us plot the points described
in Algorithm \ref{algo:main}.

\begin{figure}[H]
\begin{centering}
\includegraphics[scale=0.3]{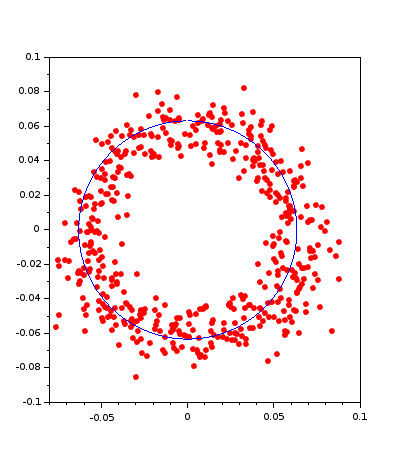}\includegraphics[scale=0.3]{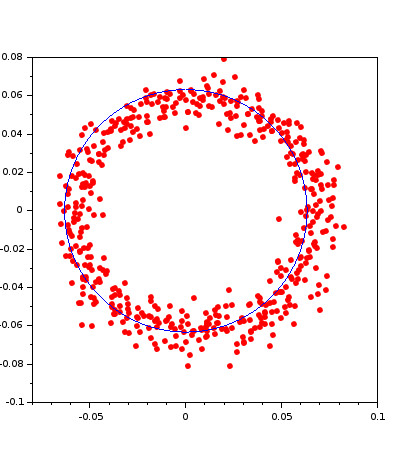}\includegraphics[scale=0.3]{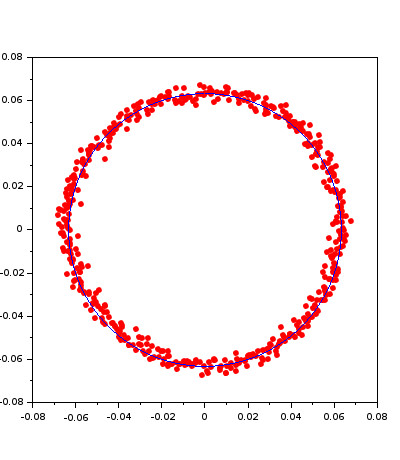}
\par\end{centering}
\caption{Points whose coordinates are entries of the eigenvectors. Probabilities:
0.3, 0.5, 0.9.}

\label{fig:eigenpoints}
\end{figure}

The circle-like shape follows indeed the layout we are looking to
reconstruct. This phenomenon can be explained in terms of angles between
spaces which appear in our proofs. Also, notice that the points that
are in the wrong position are not a major part. That can be explained
in terms of rank correlation coefficients.

A rank correlation coefficient measures the degree of similarity between
two lists, and can be used to assess the significance of the relation
between them. One can see the rank of one list as a permutation of
the rank of the other. Statisticians have used a number of different
measures of closeness for permutations. Some popular rank correlation
statistics are Kendall's $\tau$, Kendall distance, and Spearman's
footrule. There are several other metrics, and for different situations
some metrics are preferable. For a deeper discussion on metrics on
permutations we recommend \cite{Diaconis}. 

To count the total number of inversions in $\sigma$ one can use 
\[
D(\sigma)=\sum_{i<j}\textbf{1}_{\sigma(i)>\sigma(j)}\text{ (Kendall Distance)}
\]
First we define a refined version of the Kendall distance. This version
counts inverted pairs whose indices are at least $k$ positions apart.
First note that, for a permutation $\sigma$, we can rewrite $D(\sigma)$
as 
\[
D(\sigma)=|\{(i,j):\sigma(j)<\sigma(i)\mbox{ and }i<j\}|.
\]
Given a permutation $\sigma$ and an index $k\geq1$, let 
\[
D_{k}(\sigma)=|\{(i,j):\sigma(j)<\sigma(i)\textrm{ and }i+k\leq j\ \text{and }i-j\ge k\ \textrm{mod}\ n\}|.
\]
Thus, $D_{k}$ counts the number of inverted pairs where the vertices
have jumped at least $k$ positions from their original order. In
particular, $D_{1}(\sigma)=D(\sigma)$. The module in the definition
is used to access the circular structure of the graph we consider.

Consider the eigenvectors $x$ and $y$ for $\ensuremath{\lambda_{2}(M)}$
and $\lambda_{3}(M)$, respectively. As we will see, the set of points
$z_{i}=(x_{i},y_{i})$ have coordinates on a circle in $\mathbb{R}^{2}$.
Let $\varphi(x)$ be the angular coordinate for $x\in\mathbb{R}^{2}$.
A crucial observation is that $\left\{ \varphi(z_{i})\right\} _{i=1}^{n}$
is an increasing sequence. That means that the order of $\varphi(z_{i})$
provides the correct order for the vertices in the model graph. Similarly,
we can consider the eigenvectors $\widehat{x}$ and $\widehat{y}$
for $\ensuremath{\lambda_{2}(\widehat{M})}$ and $\lambda_{3}(\widehat{M})$,
respectively. Here, $\left\{ \varphi(\widehat{z_{i}})\right\} _{i=1}^{n}$
does not necessarily form an increasing sequence. Thus, we can construct
a permutation $\sigma$ of indices such that $\ensuremath{\sigma(i)>\sigma(j)}$
if and only if $\ensuremath{\varphi(\widehat{z_{i}})\geq\varphi(\widehat{z_{j}})}$.

In view of the last observations, the permutation $\sigma$ has a
neat interpretation in terms of $D_{k}$: $D_{k}(\sigma)$ counts
the pairs in $\hat{z}$ that disagree with the order induced by $z$
by at least $k$ positions. I.e, the permutation $\sigma$ of Algorithm
\ref{algo:main} has $D_{k}(\sigma)$ pairs of vertices in the wrong
order. Fortunately, the next Theorem bounds the number of such pairs.
\begin{thm}
\label{thm:firstBound}Let $\sigma$ be the permutation returned by
Algorithm \ref{algo:main} for a random circulant graph. Let $k\in\Omega(n^{\beta})$
and $|N|=cn$, for a constant $c>0$. Then it holds $D_{k}(\sigma)\in\mathcal{O}(n^{5-4\beta})$
with probability $1-n^{-3}$. 
\end{thm}

In fact, we prove a more general version of Theorem \ref{thm:firstBound}
where we allow the edge density to be variable.
\begin{thm}
\label{thm:inverted_pairs_bound}Let $\sigma$ be the permutation
returned by Algorithm \ref{algo:main} for a random circulant graph
with model satisfying $|N|=cn^{\gamma}$, for a constant $c>0$. Let
$k\in\Omega(n^{\beta}).$ Then we have $D_{k}(\sigma)\in\mathcal{O}(n^{11-6\gamma-4\beta})$
with probability $1-n^{-3}$. 
\end{thm}

Furthermore, depending on the parameters $\gamma$ and $\beta$ we
can improve the bounds of the last Theorems, as shown in the next
result. 
\begin{thm}
\label{thm:invertet_pairs_bound2}Let $\sigma$ be the permutation
returned by Algorithm \ref{algo:main} for a random circulant graph
with model satisfying $|N|=cn^{\gamma}$, for a constant $c>0$. Let
$k\in\Omega(n^{\beta}).$ Then we have $D_{k}(\sigma)\in\mathcal{O}(n^{\frac{13-6\gamma-2\beta}{3}})$
with probability $1-n^{-3}$. 
\end{thm}
Notice that the last result shows that there is a trade off between
how far vertices can jump and the total number of such incorrectly
placed vertices. That is useful for our purpose to establish metrics
on the correctness of the rank. For example, consider the worst case
of Theorem \ref{thm:invertet_pairs_bound2} when all pairs are incorrect.
Assuming $\gamma=1$, the number of pairs that drift less than $k$
positions apart is $\binom{n}{2}-D_{k}$. If we take $\beta>1/2$
in Theorem \ref{thm:invertet_pairs_bound2} , we obtain that $\binom{n}{2}-D_{k}$
is asymptotically equivalent to $n^{2}$ as $n\rightarrow\infty$.
That means almost no vertex will drift more than $n^{1/2}$ slots
from its correct position. 

Finally, the next theorem shows that the permutation returned by Algorithm
\ref{algo:main} is well behaved in terms of the usual Kendall distance.
\begin{thm}
\label{thm:error_order}Let $\sigma$ be the permutation returned
by Algorithm \ref{algo:main} for a random circulant graph with model
satisfying $|N|=cn^{\gamma}$ and $1\ge\gamma>0$. Then $D(\sigma)\in\mathcal{O}(n^{(15-6\gamma)/5})$
with probability $1-n^{-3}$. 
\end{thm}

To prove the results, our technique uses Singular Value Decomposition
and angles between subspaces, which require expressions for the eigenvalues
and eigenvectors of the model matrix. Fortunately, circulant matrices
have known spectrum and, as we will see, there is a specific pair
of eigenvectors carrying the desired information about the structure
of the graph, providing the correct label of vertices. Moreover, consecutive
entries of the eigenvectors differ significantly enough so that a
small perturbation will have limited effect on the labels. Further,
in Section \ref{sec:SVD-and-angles} we show that those eigenvectors
are close to the eigenvectors of the random graph. In Section \ref{sec:Bounds-and-proofs}
we perform the qualitative analysis of the problem proving the main
results.

\section{\label{sec:SVD-and-angles}SVD and angles between subspaces}

The definition of an angle between two vectors can be extended to
angles between subspaces.
\begin{defn}
Let $\mathcal{X}\subset\mathbb{R}^{n}$ and $\mathcal{Y}\subset\mathbb{R}^{n}$
be subspaces with $\textrm{dim}(\mathcal{X})=p$ and $\textrm{dim}(\mathcal{Y})=q$.
Let $m=\min(p,q)$. The principal angles 
\[
\Theta=\left[\theta_{1},\ldots,\theta_{m}\right]\textrm{, where }\theta_{k}\in\left[0,\pi/2\right],\textrm{ }k=1,\ldots,m,
\]
between $\mathcal{X}$ and $\mathcal{Y}$ are recursively defined
by 
\[
s_{k}=\cos(\theta_{k})=\max_{x\in\mathcal{X}}\max_{y\in\mathcal{Y}}\left|x^{T}y\right|=\left|x_{k}^{T}y_{k}\right|,
\]
subject to
\[
\left\Vert x\right\Vert =\left\Vert y\right\Vert =1,\textrm{ }x^{T}x_{i}=0,\text{}y^{T}y_{i}=0,\text{ for }i=1,\ldots,k-1.
\]
The vectors $\left\{ x_{1},\ldots,x_{m}\right\} $ and $\left\{ y_{1},\ldots,y_{m}\right\} $
are called the principal vectors for $\mathcal{X}$ and $\mathcal{Y}$. 
\end{defn}

The principal angles and principal vectors can be characterized in
terms of a Singular Value Decomposition. That provides a constructive
form for the principal vectors, which is what we use in the proofs.
That is the subject of the next Theorem proved in \cite{Golub}.
\begin{thm}
\label{thm:principalAngles}Let the columns of the matrices $X\in\mathbb{R}^{n\times p}$
and $Y\in\mathbb{R}^{n\times q}$ form an orthonormal bases for the
subspaces $\mathcal{X}$ and $\mathcal{Y}$, respectively. Consider
the singular value decomposition 
\[
X^{T}Y=U\Sigma V^{T},
\]
where $U$ and $V$ are unitary matrices and $\Sigma$ is a $p\times q$
diagonal matrix with real diagonal entries $s_{1},\ldots,s_{m}$ in
nonincreasing order with $m=\min(p,q)$. Then 
\[
\cos\Theta=\left[s_{1},\ldots,s_{m}\right],
\]
where $\Theta$ denotes the vector of principal angles between $\mathcal{X}$
and $\mathcal{Y}$. Furthermore, the principal vectors for $\mathcal{X}$
and $\mathcal{Y}$ are given by the first $m$ columns of $XU$ and
$YV$.
\end{thm}

In \cite{DK-variant}, the authors prove a variant of Davis-Kahan
Theorem, which gives an upper bound for the sine of the principal
angles between subspaces in terms of eigenvalues of the matrices whose
columns are bases for the subspaces. The original version of Davis-Kahan
\cite{DavisKahan} relies on an eigenvalue separation condition for
those matrices. However, these conditions are not necessarily met
by the eigenvalues of a random matrix. That is the reason we use a
different version of Davis-Kahan Theorem. We recast the result here
for the eigenvalues of interest of our problem. Here $\left\Vert \cdot\right\Vert _{F}$
denotes the Frobenius norm.
\begin{thm}
\label{thm:DK_Variant}Let $M,\widehat{M}\in\mathbb{R}^{n\times n}$
be symmetric matrices, with eigenvalues $\lambda_{1}\geq\ldots\geq\lambda_{n}$
and $\hat{\lambda}_{1}\geq\ldots\geq\hat{\lambda}_{n}$, respectively.
Let $\lambda_{i}$ and $\hat{\lambda_{i}}$ have corresponding unitary
eigenvectors $v_{i}$ and $\hat{v_{i}}$. Let $\min\left(\lambda_{1}-\lambda_{2},\lambda_{3}-\lambda_{4}\right)>0$,
define $V=\left[\begin{array}{cc}
v_{2} & v_{3}\end{array}\right]$ and $\widehat{V}=\left[\begin{array}{cc}
\widehat{v_{2}} & \widehat{v_{3}}\end{array}\right]$. Let $\Theta$ be a $2\times2$ diagonal matrix whose diagonal contains
the principal angles between the subspaces spanned by the columns
of $V$ and $\widehat{V}$. Then
\[
\left\Vert \sin\Theta\right\Vert _{F}\leq\frac{2\min\left(\sqrt{2}\left\Vert M-\widehat{M}\right\Vert ,\left\Vert M-\widehat{M}\right\Vert _{F}\right)}{\min\left(\lambda_{1}-\lambda_{2},\lambda_{3}-\lambda_{4}\right)}.
\]
\end{thm}

\section{\label{sec:Bounds-and-proofs}Bounds and proofs of the main results}

To prove the main theorems, we need to bound the differences $\lambda_{1}-\lambda_{2}$
and $\lambda_{3}-\lambda_{4}$. Fortunately, the spectrum of circulant
graphs is well known, see for example \cite{Hardness_results}, so
we do not need to compute it.

The four largest eigenvalues of $H$ can be expressed as follows 
\[
\lambda_{1}=\sum\limits _{k\in N}2a_{k},
\]
\[
\lambda_{2}=\lambda_{3}=\sum\limits _{k\in N}2a_{k}\cos\left(\frac{2k\pi}{n}\right),\text{ and}
\]
\[
\lambda_{4}=\sum\limits _{k\in N}2a_{k}\cos\left(\frac{4k\pi}{n}\right).
\]
Their corresponding unitary eigenvectors are 
\[
v_{1}=\frac{1}{\sqrt{n}}(1,1,\dots,1)^{T},
\]
\[
v_{2}=\frac{2}{\sqrt{2n}}(1,\cos(\frac{2\pi}{n}),\cos(2\frac{2\pi}{n}),\dots,\cos((n-1)\frac{2\pi}{n}))^{T},
\]
\[
v_{3}=\frac{2}{\sqrt{2n}}(0,\sin(\frac{2\pi}{n}),\sin(2\frac{2\pi}{n}),\dots,\sin((n-1)\frac{2\pi}{n}))^{T},\text{ and }
\]
\[
v_{4}=\frac{2}{\sqrt{2n}}(1,\cos(\frac{4\pi}{n}),\cos(2\frac{4\pi}{n}),\dots,\cos((n-1)\frac{4\pi}{n}))^{T}.
\]

Denote by $v^{i}$ is the $i-$th entry of the vector $v$. An important
observation is that the set of points with coordinates $(v_{2}^{i},v_{3}^{i})$
are on a circle in $\mathbb{R}^{2}$. Thus, these points describe
the correct structure of the graph, providing the correct label of
vertices. 

Throughout the paper $n$ is assumed to be large.
\begin{lem}
\label{lem:eig_bound} Let $H$ be a circulant graph of degree $d$
and order $n$ with eigenvalues $\lambda_{1}\ge\lambda_{2}=\lambda_{3}\ge\lambda_{4}.$
If $|N|=cn^{\gamma}$ for a constant $c>0$ and $1\ge\gamma>0$, there
is a constant $C_{1}>0$ and $C_{2}>0$ such that $\lambda_{1}-\lambda_{2}\ge C_{1}n^{3\gamma-2}$
and $\lambda_{3}-\lambda_{4}\ge C_{2}n^{3\gamma-2}.$ 
\end{lem}
\begin{proof}
We will show the lower bound for $\lambda_{1}-\lambda_{2}$ first.
Using the expression for the eigenvalues as above we have $\lambda_{1}-\lambda_{2}=\sum\limits _{k\in N}2a_{k}\left(1-\cos\left(\frac{2k\pi}{n}\right)\right).$
Note that $\cos(\theta)$ is a decreasing function in $\theta$ for
$\theta\in[0,\pi]$ and $\frac{2k\pi}{n}\le\pi$ for N and therefore
$\lambda_{1}-\lambda_{2}\ge\sum\limits _{k=1}^{|N|}2\left(1-\cos\left(\frac{2k\pi}{n}\right)\right).$
Using the Taylor series of $\cos(\theta)$ at $\theta=0$ we get $cos(\theta)\le1-\frac{\theta^{2}}{2}+\frac{\theta^{4}}{24}$,
thus 
\[
\begin{array}{lll}
\lambda_{1}-\lambda_{2} & \ge & 2\sum\limits _{k=1}^{|N|}\frac{2k^{2}\pi^{2}}{n^{2}}-\frac{2k^{4}\pi^{4}}{3n^{4}}\\
 & = & 2\left(\frac{2\pi^{2}}{n^{2}}\frac{2|N|^{3}+3|N|^{2}+|N|}{6}-\frac{2\pi^{4}}{3n^{4}}\frac{6|N|^{5}+15|N|^{4}+10|N|^{3}-|N|}{30}\right)\\
 & \ge & K_{1}\left(\frac{|N|^{3}}{n^{2}}-\frac{|N|^{5}}{n^{4}}\right),
\end{array}
\]
for a constant $K_{1}>0$. Therefore, there is a constant $C_{1}>0$
such that $\lambda_{1}-\lambda_{2}\ge C_{1}n^{3\gamma-2}.$ $\lambda_{3}-\lambda_{4}=\sum\limits _{k\in N}2a_{k}\left(\cos\left(\frac{2k\pi}{n}\right)-\cos\left(\frac{4k\pi}{n}\right)\right).$
Note, that $f(\theta)=\cos(\theta)-\cos(2\theta)$ is increasing for
$\theta\in[0,\frac{\pi}{4}]$, decreasing for $\theta\in[\frac{\pi}{4},\pi]$
and $\frac{2k\pi}{n}<\pi$ for $k\in N$. For this reason we will
split the sum above using the following partition of N, $N_{L}:=\{k\in N:\frac{2k\pi}{n}\le\frac{\pi}{4}\}$
and $N_{U}:=N-N_{L}.$ Let $\hat{k}=\max\{k\in N\}$, then using the
Taylor series for $f(\theta)$, we have 
\[
\begin{array}{lll}
\lambda_{3}-\lambda_{4} & = & \sum\limits _{k\in N_{L}}2a_{k}\left(\cos\left(\frac{2k\pi}{n}\right)-\cos\left(\frac{4k\pi}{n}\right)\right)+\sum\limits _{k\in N_{U}}2a_{k}\left(\cos\left(\frac{2k\pi}{n}\right)-\cos\left(\frac{4k\pi}{n}\right)\right)\\
 & \ge & \sum\limits _{k=1}^{|N_{L}|}2a_{k}\left(\cos\left(\frac{2k\pi}{n}\right)-\cos\left(\frac{4k\pi}{n}\right)\right)+\sum\limits _{k=\hat{k}-|N_{U}|+1}^{\hat{k}}2a_{k}\left(\cos\left(\frac{2k\pi}{n}\right)-\cos\left(\frac{4k\pi}{n}\right)\right)\\
 & \ge & 2\sum\limits _{k=1}^{|N_{L}|}\left(\frac{6\pi^{2}k^{2}}{n^{2}}-\frac{10\pi^{4}k^{4}}{n^{4}}\right)+2\sum\limits _{k=\hat{k}-|N_{U}|+1}^{\hat{k}}\left(\frac{6\pi^{2}k^{2}}{n^{2}}-\frac{10\pi^{4}k^{4}}{n^{4}}\right)\\
 & = & K_{2}\left(\frac{|N_{L}|^{3}}{n^{2}}-\frac{|N_{L}|^{5}}{n^{4}}\right)+K_{3}\left(\frac{\hat{k}^{3}}{n^{2}}-\frac{\hat{k}^{5}}{n^{4}}\right)-K_{4}\left(\frac{(\hat{k}-|N_{U}|)^{3}}{n^{2}}-\frac{(\hat{k}-|N_{U}|)^{5}}{n^{4}}\right),
\end{array}
\]
for nonnegative constants $K_{2},K_{3}$ and $K_{4}.$ Furthermore,
$\hat{k}\ge|N|$ and $\hat{k}\ge|N_{L}|$. The first inequality implies
that there is a constant $K_{5}$ such that $\hat{k}=K_{5}n^{\gamma}.$
Therefore, there is a constant $C_{2}>0$ with $\lambda_{3}-\lambda_{4}\ge C_{2}n^{3\gamma-2}$
\end{proof}
Using Lemma \ref{lem:eig_bound} we are able to prove an upper bound
for the deviations of the eigenvectors corresponding to the second
and third eigenvalues of the model matrix and the random matrix, respectively.
We will also need the following concentration inequality from \cite{concentration}.
\begin{lem}[Norm of a random matrix]
\label{lem:norm_rand_matrix} There is a costant $C>0$ such that
the following holds. Let $E$ be a symmetric matrix whose upper diagonal
entries $e_{ij}$ are independent random variables where $e_{ij}=1-p_{ij}$
or $-p_{ij}$ with probabilities $p_{ij}$ and $1-p_{ij},$ respectively,
where $0\le p_{ij}\le1.$ Let $\sigma^{2}=\max_{ij}p_{ij}(1-p_{ij}).$
If $\sigma^{2}\ge C\log n/n,$ then 
\[
\mathbf{P}(||E||\ge C\sigma n^{1/2})\le n^{-3}
\]
\end{lem}
Now we are able to prove the following theorem. 
\begin{thm}
\label{thm:up_bound} Let $M$ be the circulant graph model matrix
with constant probability $p$, variance $\sigma^{2}$, and $|N|=cn^{\gamma}$
for a constant $1\ge\gamma>0.$ Let $\hat{M}$ the random matrix following
the model matrix. Let $v_{2},v_{3}$ be unitary eigenvectors for $\lambda_{2}(M),\lambda_{3}(M)$
and $\hat{v}_{2},\hat{v}_{3}$ be unitary eigenvectors for $\lambda_{2}(\hat{M}),\lambda_{3}(\hat{M}).$
Let $x,$$y\in\textrm{Span}\left\{ v_{2},v_{3}\right\} $ and $\hat{x},$$\hat{y}\in\textrm{Span}\left\{ \hat{v}_{2},\hat{v}_{3}\right\} $
be the principal vectors for the principal angles between the spaces
$\textrm{Span}\left\{ v_{2},v_{3}\right\} $ and $\textrm{Span}\left\{ \hat{v}_{2},\hat{v}_{3}\right\} $.
Define the matrices $z=(x,y)$ and $\hat{z}=(\hat{x},\hat{y})$. Then
there is an absolute constant $C_{0}>0$ and such that 
\[
||z-\hat{z}||_{F}^{2}\le C_{0}\sigma n^{5-6\gamma}
\]
with probability at least $1-n^{-3}$.
\end{thm}
\begin{proof}
In view of Theorem \ref{thm:principalAngles}, consider the singular
value decomposition $\left[v_{2},v_{3}\right]^{T}\left[\hat{v}_{2},\hat{v}_{3}\right]=U\Sigma W^{T}.$
Let $\theta_{2}$ and $\theta_{3}$ denote the principal angles between
the spaces spanned by $\left\{ v_{2},v_{3}\right\} $ and $\left\{ \hat{v}_{2},\hat{v}_{3}\right\} $. 

Note that $\min(\lambda_{1}-\lambda_{2},\lambda_{3}-\lambda_{4})>0$,
thus we can apply Theorem \ref{thm:DK_Variant}. We have
\[
\begin{array}{lll}
||z-\hat{z}||_{F}^{2} & = & ||(x,y)-(\hat{x},\hat{y})||_{F}^{2}\\
 & = & (||x-\hat{x}||^{2}+||y-\hat{y}||^{2})\\
 & \le & 2(\sin^{2}(\theta_{2})+\sin^{2}(\theta_{3}))\\
 & = & 2||\sin\Theta(z,\hat{z})||_{F}^{2}\\
 & \le & 2\left(\frac{\min(\sqrt{2}||\hat{M}-M||,||\hat{M}-M||_{F})}{\min(\lambda_{1}-\lambda_{2},\lambda_{3}-\lambda_{4})}\right)^{2}.
\end{array}
\]

Now, we can view the adjacency matrix $\hat{M}$ as a perturbation
of $M,\ \hat{M}=M+E,$ where the entries of $E$ are $e_{ij}=1-p$
with probability $p$ and $-p$ with probability $1-p$, thus, $E$
is as in Lemma \ref{lem:norm_rand_matrix} and with probability at
least $1-n^{-3}$ we have $||E||\le C\sigma\sqrt{n}.$ Furthermore,
$\min(\sqrt{2}||\hat{M}-M||,||\hat{M}-M||_{F})\le\sqrt{2}||\hat{M}-M||_{op}=||E||$
and $\lambda_{i}(M)=p\lambda_{i}(H)$. Together with Lemma \ref{lem:eig_bound}
we get for some absolute constant $C_{0}>0$ 
\[
\begin{array}{lll}
||z-\hat{z}||^{2} & \le & C_{0}\sigma\left(\frac{\sqrt{n}}{n^{3\gamma-2}}\right)^{2}=C_{0}\sigma n^{5-6\gamma}.\end{array}
\]
That finishes the proof. 
\end{proof}

Now we will provide a lower bound for $||z-\hat{z}||_{F}$ in terms
of $D_{k}(\sigma)$ and eventually proof theorem \ref{thm:error_order}. 
\begin{lem}
\label{lem:asymptotic}Let $v$ and $w$ be the unitary eigenvectors
for $\lambda_{2}\left(M\right)$ and $\lambda_{3}\left(M\right)$,
respectively. Then
\[
\left(v_{i}-v_{i+k}\right)^{2}={n}^{-5}\left({32\,{\pi}^{4}{k}^{2}{i}^{2}+32\,{\pi}^{4}{k}^{3}i+8\,{\pi}^{4}{k}^{4}}\right)+\mathcal{O}\left(n^{-7}\right)\text{ and},
\]
\[
\left(w_{i}-w_{i+k}\right)^{2}={\left(3n\right)}^{-5}\left(12\,{\pi}^{2}{k}^{2}{n}^{2}-48\,{\pi}^{4}{k}^{2}{i}^{2}-48\,{\pi}^{4}{k}^{3}i-16\,{\pi}^{4}{k}^{4}\right)+\mathcal{O}\left(n^{-7}\right).
\]
\end{lem}
\begin{proof}
Notice that
\[
v_{i}-v_{i+k}=\frac{2}{\sqrt{2n}}(\cos\left(\frac{2\pi i}{n}\right)-\cos\left(\frac{2\pi\left(i+k\right)}{n}\right))\text{ and }
\]

\[
w_{i}-w_{i+k}=\frac{2}{\sqrt{2n}}(\sin\left(\frac{2\pi i}{n}\right)-\sin\left(\frac{2\pi\left(i+k\right)}{n}\right)).
\]
Now the expressions can be obtained from a simple asymptotic expansion
as $n\longrightarrow\infty$.
\end{proof}
\begin{thm}
\label{thm:lower_bound}Let $M$ be the circulant graph model matrix
with constant probability $p$ and variance $\sigma^{2}$, and $\hat{M}$
the random matrix following the model matrix. Let $v_{2},v_{3}$ be
unitary eigenvectors for $\lambda_{2}(M)$ and $\lambda_{3}(M)$,
and $\hat{v}_{2}$ and $\hat{v}_{3}$ be unitary eigenvectors for
$\lambda_{2}(\hat{M})$ and $\lambda_{3}(\hat{M}).$ Let $x,$$y\in\textrm{Span}\left\{ v_{2},v_{3}\right\} $
and $\hat{x},$$\hat{y}\in\textrm{Span}\left\{ \hat{v}_{2},\hat{v}_{3}\right\} $
be the principal vectors for the principal angles between the spaces
$\textrm{Span}\left\{ v_{2},v_{3}\right\} $ and $\textrm{Span}\left\{ \hat{v}_{2},\hat{v}_{3}\right\} $.
Define the matrices $z=(x,y)$ and $\hat{z}=(\hat{x},\hat{y})$. Then
there are constants $C_{1}>0$ and $\beta$ such that 
\[
\left\Vert z-\hat{z}\right\Vert ^{2}>C_{0}\left|R\right|\frac{n^{4\beta}}{n^{6}},
\]
 where $R=\{(i,j):\varphi\left(\hat{z_{j}}\right)\ge\varphi\left(\hat{z_{i}}\right),i+k\leq j$
and $i-j\ge k\ \textrm{mod}\ n\}.$
\end{thm}
\begin{proof}
As in Theorem \ref{thm:principalAngles}, let $U\Sigma W^{T}$ be
the singular value decomposition for the matrix $\left[v_{2},v_{3}\right]^{T}\left[\hat{v}_{2},\hat{v}_{3}\right].$
Thus, $z=(v_{2},v_{3})U$ and $\hat{z}=(\hat{v}_{2},\hat{v}_{3})W$.
Let $\varphi(z_{i})$ be the angular coordinate of the point $z_{i}=(x_{i},y_{i})$.
Thus, $\left\{ \varphi(z_{i})\right\} _{i=1}^{n}$ is an increasing
sequence. Fix $k=k(n)=C(n^{\beta})$ and let 
\[
R=\{(i,j):\varphi\left(\hat{z_{j}}\right)\le\varphi\left(\hat{z_{i}}\right)\mbox{ and }i+k\leq j\text{ and }i-j\ge k\ mod\ n\}.
\]
Then $R$ is the set of pairs in $\hat{z}$ that disagree with the
order induced by $z$ by at least $k$ positions in both directions
on the cycle. Now we can write 
\[
2n\Vert z-\hat{z}\Vert^{2}=\sum_{i=1}^{n}\sum_{j=1}^{n}\Vert z_{i}-\hat{z_{i}}\Vert^{2}+\Vert z_{j}-\hat{z}_{j}\Vert^{2}\geq\sum_{(i,j)\in R}\Vert z_{i}-\hat{z_{i}}\Vert^{2}+\Vert z_{j}-\hat{z}_{j}\Vert^{2}.
\]

Since $\varphi\left(\hat{z_{j}}\right)\le\varphi\left(\hat{z_{i}}\right)$
and $\varphi\left(z_{j}\right)>\varphi\left(z_{i}\right)$, the minimum
contribution of each term in the sum happens in the median point 
\[
\hat{z_{i}}=\hat{z_{j}}=\frac{z_{i}+z_{j}}{2}.
\]
Thus, we have

\begin{eqnarray*}
2n\left\Vert z-\hat{z}\right\Vert ^{2} & > & \sum_{(i,j)\in R}\Vert z_{i}-\frac{z_{i}+z_{j}}{2}\Vert^{2}+\Vert z_{j}-\frac{z_{i}+z_{j}}{2}\Vert^{2}\\
 & = & \sum_{(i,j)\in R}\frac{\Vert z_{i}-z_{j}\Vert^{2}}{2}
\end{eqnarray*}

Now set $v=v_{2}$ and $w=v_{3}$ and notice that $(z_{i}-z_{j})U^{T}=(v_{i}-v_{j},w_{i}-w_{j})$
and, by the definition of $R$, $||z_{i}-z_{j}||^{2}$ is minimum
for $j=i+k.$ Thus we write

\begin{eqnarray*}
2n\left\Vert z-\hat{z}\right\Vert ^{2} & > & \sum_{(i,j)\in R}\frac{\Vert z_{i}-z_{i+k}\Vert^{2}}{2}\\
 & = & \sum_{(i,j)\in R}\frac{||(z_{i}-z_{i+k})U^{T}||^{2}}{2}\\
 & = & \sum_{(i,j)\in R}\frac{(v_{i}-v_{i+k})^{2}+(w_{i}-w_{i+k})^{2}}{2}.
\end{eqnarray*}
By Lemma \ref{lem:asymptotic}, for $n$ large enough we obtain

\[
2n\left\Vert z-\hat{z}\right\Vert ^{2}>C\sum_{(i,j)\in R}\frac{k^{4}}{n^{5}},
\]
for an absolute constant $C>0$. 

Now $k=k(n)=\Omega(n^{\beta})$, so there exists a constant $c$ such
that for $n$ large enough $k\geq cn^{\beta}$. We can bound

\[
\left\Vert z-\hat{z}\right\Vert ^{2}>C_{0}\left|R\right|\frac{n^{4\beta}}{n^{6}},
\]
for a constant $C_{0}>0$.
\end{proof}

Now, the proof of Theorem \ref{thm:inverted_pairs_bound} easily follows
from Theorem \ref{thm:lower_bound} and Theorem \ref{thm:up_bound}.
First, we make an observation about the order given by Algorithm \ref{algo:main}.
Let $v_{2},v_{3}$ be unitary eigenvectors for $\lambda_{2}(M)$ and
$\lambda_{3}(M)$, and $\hat{v}_{2}$ and $\hat{v}_{3}$ be unitary
eigenvectors for $\lambda_{2}(\hat{M})$ and $\lambda_{3}(\hat{M}).$
Let $U\Sigma W^{T}$ be the singular value decomposition for the matrix
$\left[v_{2},v_{3}\right]^{T}\left[\hat{v}_{2},\hat{v}_{3}\right].$
Let $\widehat{x}=\hat{v}_{2}$ and $\widehat{y}=\hat{v}_{3}$ as in
Algorithm \ref{algo:main}, and let $\ensuremath{\varphi\left((\widehat{x_{i}},\widehat{y_{i}})\right)}$
be the angular coordinate of the point $\ensuremath{(\widehat{x_{i}},\widehat{y_{i}})}$.
Define the matrices $z=(v_{2}\ v_{3})U$ and $\hat{z}=(\hat{v}_{2}\ \hat{v}_{3})W$.
Finally, let 
\[
R=\{(i,j):\varphi\left(\hat{z_{j}}\right)\le\varphi\left(\hat{z_{i}}\right)\mbox{ and }i+k\leq j\text{ and }i-j\ge k\ \textrm{mod}\ n\}.
\]

Notice that since $W$ is a rotation matrix, it holds
\[
\ensuremath{\varphi\left((\widehat{x_{i}},\widehat{y_{i}})\right)-\ensuremath{\varphi\left((\widehat{x_{j}},\widehat{y_{j}})\right)}=\varphi\left(\hat{z_{i}}\right)-\ensuremath{\varphi\left(\hat{z_{j}}\right).}}
\]
Thus the order induced by the row vectors of $\hat{z}$ is the same
as the order induced by the row vectors of $\left[\hat{v}_{2},\hat{v}_{3}\right]$.
That implies $D_{k}(\sigma)=\left|R\right|$, where $\sigma$ is the
permutation returned by the Algorithm \ref{algo:main}. Therefore,
we proceed bounding $\left|R\right|$.

\begin{proof}
(Theorem \ref{thm:inverted_pairs_bound}) By theorems \ref{thm:lower_bound}
and \ref{thm:up_bound} we have
\[
C_{0}|R|\frac{n^{4\beta}}{n^{6}}<||z-\hat{z}||^{2}\le\bar{C_{0}}n^{5-6\gamma},
\]
where $C_{0}$ and $\bar{C_{0}}$ are positive constants and the upper
bound holds with probability at least $1-n^{-3}.$ Therefore, there
is a constant $C>0$ such that
\[
|R|<Cn^{11-6\gamma-4\beta},
\]
with probability at least $1-n^{-3}.$
\end{proof}

The proof of Theorem \ref{thm:invertet_pairs_bound2} is similar to
the last one but uses a different trick to get another lower bound. 
\begin{proof}
(Theorem \ref{thm:invertet_pairs_bound2}) As in Theorem \ref{thm:principalAngles},
let $U\Sigma W^{T}$ be the singular value decomposition for the matrix
$\left[v_{2},v_{3}\right]^{T}\left[\hat{v}_{2},\hat{v}_{3}\right].$
Thus, $z=(v_{2},v_{3})U$ and $\hat{z}=(\hat{v}_{2},\hat{v}_{3})W$.
Let $\varphi(z_{i})$ be the angular coordinate of the point $z_{i}=(x_{i},y_{i})$.
Thus, $\left\{ \varphi(z_{i})\right\} _{i=1}^{n}$ is an increasing
sequence. Fix $k=k(n)=C(n^{\beta})$ and let 
\[
R=\{(i,j):\varphi\left(\hat{z_{j}}\right)\le\varphi\left(\hat{z_{i}}\right)\mbox{ and }i+k\leq j\text{ and }i-j\ge k\ mod\ n\}.
\]
Then $R$ is the set of pairs in $\hat{z}$ that disagree with the
order induced by $z$ by at least $k$ positions in both directions
on the cycle. Now we can write 
\[
2n\Vert z-\hat{z}\Vert^{2}=\sum_{i=1}^{n}\sum_{j=1}^{n}\Vert z_{i}-\hat{z_{i}}\Vert^{2}+\Vert z_{j}-\hat{z}_{j}\Vert^{2}\geq\sum_{(i,j)\in R}\Vert z_{i}-\hat{z_{i}}\Vert^{2}+\Vert z_{j}-\hat{z}_{j}\Vert^{2}.
\]

Since $\varphi\left(\hat{z_{j}}\right)\le\varphi\left(\hat{z_{i}}\right)$
and $\varphi\left(z_{j}\right)>\varphi\left(z_{i}\right)$, the minimum
contribution of each term in the sum happens in the median point 
\[
\hat{z_{i}}=\hat{z_{j}}=\frac{z_{i}+z_{j}}{2}.
\]
Thus, we have

\begin{eqnarray}
2n\left\Vert z-\hat{z}\right\Vert ^{2} & > & \sum_{(i,j)\in R}\Vert z_{i}-\frac{z_{i}+z_{j}}{2}\Vert^{2}+\Vert z_{j}-\frac{z_{i}+z_{j}}{2}\Vert^{2}\nonumber \\
 & = & \sum_{(i,j)\in R}\frac{\Vert z_{i}-z_{j}\Vert^{2}}{2}\label{ineq_lower}
\end{eqnarray}
Now, let $n_{i}$ denote the number of pairs $(i,j)$ and label them
$(i,j_{i_{1}}),\dots,(i,j_{i_{n_{i}}})$ for $i=1,\dots,n$ and therefore
$\sum_{i=1}^{n}n_{i}=|R|$. Furthermore, for a fixed i we obtain the
minimum $||z_{i}-z_{j}||^{2}$whenever $j-i$ is minimum. By definition
of $R$ the sum
\[
\sum_{t=1}^{n_{i}}\frac{||z_{i}-z_{j}||^{2}}{2}
\]
is minimum whenever $j_{i_{1}}=i+k,\ j_{i_{2}}=i+k+1,\dots,j_{i_{n_{i}}}=i+k+n_{i}-1.$
Therefore, setting $v=v_{2}$ and $w=v_{3}$ , inequality \ref{ineq_lower}
becomes
\[
2n||z-\hat{z}||^{2}>\sum_{i=1}^{n}\sum_{t=0}^{n_{i}-1}\frac{||z_{i}-z_{i+k+t}||^{2}}{2}=\sum_{i=1}^{n}\sum_{t=0}^{n_{i}-1}\frac{(v_{i}-v_{i+k+t})^{2}+(w_{i}-w_{i+k+t})^{2}}{2}.
\]
By Lemma \ref{lem:asymptotic}, for n large enough,
\[
2n||z-\hat{z}||^{2}>C_{1}\sum_{i=1}^{n}\sum_{t=0}^{n_{i}-1}\frac{(k+t)^{4}}{n^{5}}>C_{1}\frac{k^{2}}{n^{5}}\sum_{i=1}^{n}\sum_{t=0}^{n_{i}-1}t^{2}
\]
for a constant $C_{1}>0.$ Therefore, there is a constant $C_{2}>0$
such that $2n||z-\hat{z}||^{2}>C_{2}\frac{k^{2}}{n^{5}}\sum_{i=1}^{n}n_{i}^{3}.$
Now, recall that two $p$-norms are related by $||x||_{p}\le n^{\frac{1}{p}-\frac{1}{q}}||x||_{q}.$
Taking $p=1$ and $q=3$, we obtain
\[
\sum_{i=1}^{n}n_{i}\le n^{\frac{2}{3}}(\sum_{i=1}^{n}n_{i}^{3})^{\frac{1}{3}}
\]
which allows us to rewrite inequality \ref{ineq_lower} as
\[
2n||z-\hat{z}||>C_{1}k^{2}|R|^{3}n^{-7}.
\]
Combining this inequality with the upper bound of Theorem \ref{thm:up_bound}
and using $k=k(n)=\Omega(n^{\beta})$, we obtain a constant $C_{2}>0$
such that, with probability $1-n^{-3},$
\[
|R|<C_{2}n^{\frac{13-6\gamma-2\beta}{3}}
\]
and therefore $D_{k}(\sigma)\in\mathcal{O}(n^{\frac{13-6\gamma-2\beta}{3}})$.
\end{proof}
Eventually, we give a proof for Theorem \ref{thm:error_order}.
\begin{proof}
(Theorem \ref{thm:error_order}) Fix $k=n^{(10-6\gamma)/5}$ and define
\[
R=\{(i,j):\varphi\left(\hat{z_{j}}\right)\le\varphi\left(\hat{z_{i}}\right)\mbox{ and }i+k\leq j\text{ and }i-j\ge k\ \textrm{mod}\ n\}
\]

and
\[
R^{C}=\{(i,j):\varphi\left(\hat{z_{j}}\right)\le\varphi\left(\hat{z_{i}}\right)\mbox{ and }j<i+k\text{ or }i-j<k\ \textrm{mod}\ n)\}.
\]
By Theorem \ref{thm:inverted_pairs_bound}, taking $\beta=\frac{10-6\gamma}{5},$
there is a constant $C>0$ so that, for large enough $n$,
\[
|R|=D_{k}(\sigma)\le Cn^{11-6\gamma-4\beta}=Cn^{(15-6\gamma)/5}.
\]
Furthermore, for each index $i$ there are at most $2k$ pairs $(i,j)$
in $R^{C},$thus
\[
|R^{C}|\le2kn=2n^{\beta+1}=2n^{(15-6\gamma)/5}
\]
and therefore $D(\sigma)=|R|+|R^{C}|\le(C+2)n^{(15-6\gamma)/5}$,
as required.
\end{proof}

\end{document}